\documentclass[11pt]{article}
\usepackage{amsfonts}
\usepackage{amssymb}
\usepackage{epsfig}
\usepackage{graphicx}
\usepackage{float}
\newtheorem{theorem}{Theorem}[section]
\newtheorem{proposition}[theorem]{Proposition}

\newtheorem{definition}{Definition}[section]
\newtheorem{remark}{Remark}[section]
\newtheorem{example}{Example}[section]
 \usepackage{amssymb}
\usepackage{epsfig}
\usepackage{amsmath}
\usepackage{amsfonts}

\newcommand{\R}{{{\Bbb R}}}

\def\qed{\hbox to 0pt{}\hfill$\rlap{$\sqcap$}\sqcup$\medbreak}

\title{Existence of minimal and maximal solutions to first--order differential equations with state--dependent deviated arguments}
\author{Rub\'en Figueroa and Rodrigo L\'opez Pouso}
\begin{document}
\maketitle

\vspace{1cm}
\begin{center}
{\large Departamento de
An\'alise Matem\'atica,\\
Facultade de Matem\'aticas,\\Universidade de Santiago de Compostela, Campus Sur,\\
15782 Santiago de
Compostela, Spain.
%\\
%Phone: 34 981 56 31 00,
%Ext. 13213 \\ FAX: 34 981 59 70 54
}
\end{center}

\begin{abstract}
We prove some new results on existence of solutions to first--order ordinary differential equations with deviating arguments. Delay differential equations are included in our general framework, which even allows deviations to depend on the unknown solutions. Our existence results lean on new definitions of lower and upper solutions introduced in this paper, and we show with an example that similar results with the classical definitions are false. We also introduce an example showing that the problems considered need not have the least (or the greatest) solution between given lower and upper solutions, but we can prove that they do have minimal and maximal solutions in the usual set--theoretic sense. Sufficient conditions for the existence of lower and upper solutions, with some examples of application, are provided too.
\end{abstract}

\section{Introduction}
Let $I_0=[t_0,t_0+L]$ be a closed interval, $r \ge 0$, and put $I_-=[t_0-r,t_0]$ and $I=I_- \cup I_0$. In this paper we are concerned with the existence of solutions for the following problem with deviated arguments:

\begin{equation}\label{p1}
\left\{
\begin{array}{ll}
x'(t)=f(t,x(t),x(\tau(t,x))) \ \mbox{for almost all (a.a.) } t \in I_0, \\
\\
x(t)=\Lambda(x)+k(t) \ \mbox{for all $t \in I_-$,}
\end{array}
\right.
\end{equation}
where $f:  I \times \R^2 \longrightarrow \R$ and $\tau:   I_0 \times \mathcal{C}(I) \longrightarrow I$ are Carath\'eodory functions, $\Lambda:  {\cal C}(I) \longrightarrow \R$ is  a continuous nonlinear operator and $k \in \mathcal{C}(I_-)$. Here ${\cal C}(J)$ denotes the set of real functions which are continuous on the interval $J$. \\

We define a solution of problem (\ref{p1}) to be a function $x \in \mathcal{C}(I)$ such that $x_{|I_0} \in AC(I_0)$ (i.e., $x_{|I_0}$ is absolutely continuous on $I_0$) and $x$ fulfills (\ref{p1}).\\

In the space $\mathcal{C}(I)$ we consider the usual pointwise partial
ordering, i.e., for $\gamma_1,\gamma_2 \in \mathcal{C}(I)$ we define $\gamma_1 \le \gamma_2$ if and only if
$\gamma_1(t) \le \gamma_2(t)$ for all $t \in I.$ A solution of (\ref{p1}), $x_*$, is a {\it minimal} (respectively, {\it maximal}) solution of (\ref{p1}) in a certain subset $Y \subset {\cal C}(I)$ if $x_* \in Y$ and the inequality $x  \le x_*$ (respectively, $x \ge x_*$) implies $x =x_*$ whenever $x$ is a solution to (\ref{p1}) and $x \in Y$. We say that $x_*$ is {\it the least} (respectively, {\it the greatest}) solution of (\ref{p1}) in $Y$ if $x_* \le x$ (respectively, $x_* \ge x$) for any other solution $x \in Y$. Notice that the least solution in a subset $Y$ is a minimal solution in $Y$, but the converse is false in general, and an analgous remark is true for maximal and greatest solutions.  \\

 Interestingly, we will show that problem (\ref{p1}) may have minimal (ma\-xi\-mal) solutions between given lower and upper solutions and not have the least (greatest) solution. This seems to be a peculiar feature of equations with deviating arguments, see \cite{figpou2} for an example with a second--order equation.  Therefore, we are obliged to distinguish between the concepts of mi\-nimal solution and least solution (or maximal and greatest solutions), unfortunately identified in the literature on lower and upper solutions.  \\

First--order differential equations with state--dependent deviated arguments have received a lot of attention in the last years. We can cite the recent papers \cite{arino}, \cite{dyk}, \cite{dykjan}, \cite{jan}, \cite{jan1}, \cite{walther} which deal with existence results for this kind of problems.
For the qualitative study of this type of problems we can cite the survey of Hartung {\it et al.} \cite{har} and references therein. \\
As main improvements in this paper with regard to previous works in the literature we can cite the following:

\begin{enumerate}
\item[$(1)$] The deviating argument $\tau$ depends at each moment $t$ on the global behaviour of the solution, and not only on the values that it takes at the instant $t$.
\item[$(2)$] Delay problems, which correspond to differential equations of the form $x'(t)=f(t,x(t),x(t-r))$ along with a functional start condition, are included in the framework of problem (\ref{p1}).
This is not allowed in papers \cite{dyk}, \cite{dykjan}, \cite{jan} and \cite{jan1}.
\item[$(3)$] No monotonicity conditions are required for the functions $f$ and $\tau$, and they need not be continuous with respect to their first variable.
\end{enumerate}

This paper is organized as follows. In Section $2$ we state and prove the main results in this paper, which are two existence results for problem (\ref{p1}) between given lower and upper solutions.  The first result ensures the existence of maximal and minimal solutions, and the second one establishes the existence of the greatest and the least solutions in a particular case. The concepts of lower and upper solutions introduced in Section 2 are new, and we show with an example that our existence results are false if we consider lower and upper solutions in the usual sense. We also show with an example that our problems need not have the least or the greatest solution between given lower and upper solutions. In Section $3$ we prove some results on the existence of lower and upper solutions with some examples of application.

\section{Main results}

\noindent We begin this section by introducing adequate new definitions of lower and upper solutions for problem (\ref{p1}).

Notice first that $\tau(t,\gamma) \in I=I_-\cup I_0$ for all $(t,\gamma) \in I_0 \times \mathcal{C}(I)$, so for each $t \in I_0$ we can define
$$
\tau_*(t)=\inf_{\gamma \in \mathcal{C}(I)} \tau(t,\gamma) \in I, \quad \tau^*(t)=\sup_{\gamma \in \mathcal{C}(I)} \tau(t,\gamma) \in I.$$

\begin{definition}\label{sub-sobre} We say that $\alpha,\beta \in \mathcal{C}(I)$, with $\alpha \le \beta$ on $I$, are a lower and an upper solution for problem (\ref{p1}) if $\alpha_{|I_0},\beta_{|I_0} \in AC(I_0)$ and the following inequalities hold:
\begin{eqnarray}
\label{sub1}
\alpha'(t) \le \min_{\xi \in E(t)} f(t,\alpha(t),\xi) \ \mbox{for a.a. } t \in I_0, \quad \alpha  \le \inf_{\gamma \in [\alpha,\beta]} \Lambda (\gamma) +k  \ \mbox{on $I_-$,} \\
\label{sobre1}
\beta'(t) \ge \max_{\xi \in E(t)} f(t,\beta(t),\xi) \ \mbox{for a.a. } t \in I_0, \quad \beta  \ge  \sup_{\gamma \in [\alpha,\beta]} \Lambda (\gamma) + k \ \mbox{on $I_-$,}
\end{eqnarray}
where
$$E(t)=\left[\min_{s \in [\tau_*(t),\tau^*(t)]} \alpha(s),\max_{s \in [\tau_*(t),\tau^*(t)]} \beta(s)\right] \quad (t \in I_0),$$
and $[\alpha,\beta]=\{ \gamma \in {\cal C}(I) \, : \, \mbox{$\alpha \le \gamma \le \beta$ on $I$}\}.$

\end{definition}

\begin{remark}
Definition  \ref{sub-sobre} requires implicitly that $\Lambda$ be bounded in $[\alpha,\beta]$.

On the other hand, the values
$$\min_{\xi \in E(t)} f(t,\alpha(t),\xi) \quad \mbox{and} \quad \max_{\xi \in E(t)} f(t,\beta(t),\xi),$$
are really attained for almost every fixed $t \in I_0$ thanks to the continuity of $f(t,\alpha(t),\cdot)$ and $f(t,\beta(t),\cdot)$ on the compact set $E(t)$.
\end{remark}

Now we introduce the main result in this paper.

\begin{theorem}\label{main1} Assume that the following conditions hold:
\begin{enumerate}
\item[$(H_1)$] {\it (Lower and upper solutions)} There exist $\alpha,\beta \in \mathcal{C}(I)$, with $\alpha \le \beta$ on $I$, which are
 a lower and an upper solution for problem (\ref{p1}).
\item[$(H_2)$] {\it (Carath\'eodory conditions)}

\hspace{.5cm} $(H_2)-(a)$ For all $x,y \in [\min_{t \in I} \alpha(t),\max_{t \in I} \beta(t)]$ the function $f(\cdot,x,y)$ is measurable and for a.a. $t \in I_0$, all $x \in [\alpha(t),\beta(t)]$
and all $y \in E(t)$ (as defined in Definition \ref{sub-sobre}) the functions $f(t,\cdot,y)$ and $f(t,x,\cdot)$ are continuous. \\

\hspace{.5cm} $(H_2)-(b)$ For all $\gamma \in [\alpha,\beta]=\{\xi \in \mathcal{C}(I) \, : \, \alpha \le \xi \le \beta \ \mbox{on $I$}\}$ the function $\tau(\cdot,\gamma)$
is measurable and for a.a. $t \in I_0$ the operator $\tau(t,\cdot)$ is continuous in ${\cal C}(I)$ (equipped with it usual topology of uniform convergence). \\

\hspace{.5cm} $(H_2)-(c)$ The nonlinear operator $\Lambda:{\cal C}(I)\longrightarrow \R$ is continuous.

\item[$(H_3)$] {\it ($L^1-$bound)} There exists $\psi \in L^1(I_0)$ such that for a.a. $t \in I_0$, all $x \in [\alpha(t),\beta(t)]$ and all $y \in E(t)$ we have
$$
|f(t,x,y)| \le \psi(t).$$
\end{enumerate}
Then problem (\ref{p1}) has maximal and minimal solutions in $[\alpha,\beta]$.
\end{theorem}

\noindent {\bf Proof.} As usual, we consider the function
$$
p(t,x)=\left\{
\begin{array}{cl}
\alpha(t), \ &\mbox{if } x < \alpha(t),\\
x, \ &\mbox{if } \alpha(t) \le x \le \beta(t),\\
\beta(t), \ &\mbox{if } x > \beta(t),
\end{array}
\right.
$$
and the modified problem
\begin{equation}\label{paux}
\left\{
\begin{array}{ll}
x'(t)=f(t,p(t,x(t)),p(\tau(t,x),x(\tau(t,x)))) \ \mbox{for a.a. } t \in I_0, \\
\\
x(t)=\Lambda(p(\cdot,x(\cdot)))+k(t) \ \mbox{for all $t \in I_-.$}
\end{array}
\right.
\end{equation}

\noindent
{\it Claim $1$: Problem (\ref{paux}) has a nonempty and compact set of solutions.} Consider the operator $T: \mathcal{C}(I) \longrightarrow \mathcal{C}(I)$ which maps each $\gamma \in \mathcal{C}(I)$ to a continuous function $T\gamma$ defined for each $t \in I_-$ as
$$
T\gamma(t)=
\Lambda(p(\cdot,\gamma(\cdot)))+k(t),$$ and for each $t \in I_0$ as
$$T\gamma(t)=
\Lambda(p(\cdot,\gamma(\cdot)))+k(t_0)+ \int_{t_0}^t f(s,p(t,\gamma(s)),p(\tau(s,\gamma),\gamma(\tau(s,\gamma))))ds.$$
It is an elementary matter to check that $T$ is a completely con\-ti\-nuous ope\-rator from $\mathcal{C}(I)$ into itself, and therefore Schauder's Theorem ensures that $T$ has a nonempty and compact set of fixed points in $\mathcal{C}(I)$, which are exactly the solutions of problem (\ref{paux}). \\

\noindent
{\it Claim $2$: Every solution $x$ of (\ref{paux}) satisfies $\alpha \le x \le \beta$ on $I$ and, therefore, it is a solution of (\ref{p1}) in $[\alpha,\beta]$.} First, notice that if $x$ is a solution of (\ref{paux}) then $ p(\cdot,x(\cdot))
\in[\alpha,\beta]$. Hence the definition of lower solution implies that for all $t \in I_-$ we have
$$\alpha(t) \le \Lambda(p(\cdot,x(\cdot)))+k(t)=x(t).$$

Assume now, reasoning by contradiction, that $x \ngeq \alpha$ on $I_0$. Then we can find $\hat{t}_0?\in [t_0,t_0+L)$ and $\varepsilon>0$ such that $\alpha(\hat{t}_0)=x(\hat{t}_0)$ and
\begin{equation}
\label{cont}
\alpha(t) >x(t) \, \, \mbox{for all $t \in [\hat{t}_0,\hat{t}_0+\varepsilon].$}
\end{equation}
Therefore, for all $t \in [\hat{t}_0,\hat{t}_0+\varepsilon]$ we have $p(t,x(t))=\alpha(t)$ and
$$p(\tau(t,x),x(\tau(t,x))) \in [\alpha(\tau(t,x)),\beta(\tau(t,x))] \subset E(t),$$ so for a.a. $s \in [\hat{t}_0,\hat{t}_0+\varepsilon]$ we have
$$
\alpha'(s) \le f(s,p(s,x(s)),p(\tau(s,x),x(\tau(s,x)))).$$
Hence for $t \in [\hat{t}_0,\hat{t}_0+\varepsilon]$ we have
$$ \alpha(t)-x(t)=\int_{\hat{t}_0}^{t} \alpha'(s) \, ds - \int_{\hat{t}_0}^{t} f(s,p(s,x(s)),p(\tau(s,x),x(\tau(s,x)))) \, ds \le 0,$$
a contradiction with (\ref{cont}).

 Similar arguments prove that all solutions $x$ of (\ref{paux}) obey $x \le \beta$ on $I$.
\noindent
{\it Claim $3$: The set of solutions of problem (\ref{p1}) in $[\alpha,\beta]$ has maximal and minimal elements.} The set
$$
\mathcal{S}=\{x \in \mathcal{C}(I) \, : \, \mbox{$x$ is a solution of (\ref{p1}) }, \ x \in [\alpha,\beta]\}$$
is nonempty and compact in $\mathcal{C}(I)$, beacuse it coincides with the set of fixed points of the operator $T$. Then, the real--valued continuous ma\-pping
$$
x \in \mathcal{S} \longmapsto \mathcal{I}(x)=\int_{t_0}^{t_0+L} x(s) \, ds$$
attains its maximum and its minimum, that is, there exist $x^*,x_* \in \mathcal{S}$ such that
\begin{equation}\label{int}
\mathcal{I}(x^*)=\max\{\mathcal{I}(x)\, : \, x \in \mathcal{S}\}, \quad \mathcal{I}(x_*)=\min\{\mathcal{I}(x)\, : \, x \in \mathcal{S}\}.\end{equation}
Now, if $x \in \mathcal{S}$ is such that $x \ge x^*$ on $I$ then we have $\mathcal{I}(x) \ge \mathcal{I}(x^*)$ and, by (\ref{int}), $\mathcal{I}(x) \le \mathcal{I}(x^*)$. So we conclude that $\mathcal{I}(x) = \mathcal{I}(x^*)$ which, along with $x \ge x^*$, implies that $x=x^*$ on $I$. Hence $x^*$ is a maximal element of $\mathcal{S}$. In the same way we can prove that $x_*$ is a minimal element.
\qed

One might be tempted to follow the standard ideas with lower and upper solutions to define a lower solution of (\ref{p1}) as some function $\alpha$ such that
\begin{equation}
\label{sub2}
\alpha'(t) \le f(t,\alpha(t),\alpha(\tau(t,\alpha))) \ \mbox{for a.a. } t \in I_0,
\end{equation}
and an upper solution as some function $\beta$ such that
\begin{equation}
\label{sobre2}
\beta'(t) \ge f(t,\beta(t),\beta(\tau(t,\beta))) \ \mbox{for a.a. } t \in I_0.
\end{equation}

These definitions are not adequate to ensure the existence of solutions of (\ref{p1}) between given lower and upper solutions, as we show in the following example.

\begin{example}\label{usual} Consider the problem with delay
\begin{equation}\label{exusual}
x'(t)=-x(t-1), \ \mbox{for a.a. $t \in [0,1]$}, \quad x(t)=k(t)=-t, \ \mbox{for $t \in [-1,0]$}.
\end{equation}

Notice that functions $\alpha(t)=0$ and $\beta(t)=1$, $t\in [-1,1]$, are lower and upper solutions in the usual sense for problem (\ref{exusual}). However, if $x$ is a solution for problem (\ref{exusual}) then for a.a. $t \in [0,1]$ we have $$
x'(t)=-x(t-1)=-k(t-1)=-[-(t-1)]=t-1,$$
so for all $t \in [0,1]$ we compute
$$x(t)=x(0)+\int_0^t (s-1) \, ds = \dfrac{t^2}{2}-t,$$
and then $x(t) < \alpha(t)$ for all $t \in (0,1]$. Hence (\ref{exusual}) has no solution at all between $\alpha$ and $\beta$.
\end{example}

\begin{remark}\label{rem1} Notice that inequalities (\ref{sub1}) and (\ref{sobre1}) imply (\ref{sub2}) and (\ref{sobre2}),
so lower and upper solutions in the sense of Definition \ref{sub-sobre} are lower and upper solutions in the usual sense, but the converse is false in general.

Definition \ref{sub-sobre} is probably the best possible for (\ref{p1}) because it reduces to some definitions that one can find in the literature in connection with particular cases of (\ref{p1}). Indeed, when the function $\tau$ does not depend on the second variable then for all $t \in I_0$ we have $E(t)=[\alpha(\tau(t)),\beta(\tau(t))]$ in Definition \ref{sub-sobre}. Therefore, if $f$ is nondecreasing with respect to its third variable, then Definition \ref{sub-sobre} and the usual definition of lower and upper solutions are the same (we will use this fact in the proof of Theorem \ref{main2}). If, in turn, $f$ is nonincreasing with respect to its third variable, then Definition \ref{sub-sobre} coincides with the usual definition of coupled lower and upper solutions (see for example \cite{jan}).
\end{remark}

In general, in the conditions of Theorem \ref{main1} we cannot expect problem (\ref{p1}) to have the extremal solutions  in $[\alpha,\beta]$ (that is, the greatest and the least solutions in $[\alpha,\beta]$). This is justified by the following example.

\begin{example}\label{extremal1} Consider the problem
\begin{equation}\label{ex1}
x'(t)=f(t,x(t),x(\tau(t))) \ \mbox{for a.a. } t \in I_0=\left[-\frac{\pi}{2},\pi\right], \quad x\left(-\frac{\pi}{2}\right)=0,
\end{equation}
where
$$
f(t,x,y)=\left\{
\begin{array}{cl}
1, \ &\mbox{if } y < -1, \\
-y, \ &\mbox{if } -1 \le y \le 1, \\
-1, \ &\mbox{if } y >1,
\end{array}
\right.
$$
and $\tau(t)= \frac{\pi}{2}-t$. \\

First we check that $\alpha(t)=-t-\frac{\pi}{2}=-\beta(t), \ t \in I_0,$ are lower and upper solutions for problem (\ref{ex1}). The definition of $f$ implies that for all $(t,x,y) \in I_0\times \R^2$ we have $|f(t,x,y)|\le 1,$
so for all $t \in I_0$ we have
$$
\min_{\xi \in E(t)} f(t,\alpha(t),\xi)\ge -1 = \alpha'(t) \, \, \mbox{and} \, \, \max_{\xi \in E(t)} f(t,\beta(t),\xi)\le 1 = \beta'(t),$$
where, according to Definition \ref{sub-sobre},$$
E(t)=\left[\alpha\left(\frac{\pi}{2}-t\right),\beta\left(\frac{\pi}{2}-t\right)\right]=
\left[t-\pi,\pi-t\right].$$

Moreover, $\alpha(-\frac{\pi}{2})=\beta(-\frac{\pi}{2})=0$, so $\alpha$ and $\beta$ are, respectively, a lower and an upper solution for (\ref{ex1}), and then condition $(H_1)$ of Theorem \ref{main1}is fulfilled. As conditions $(H_2)$ and $(H_3)$ are also satisfied (take, for example, $\psi \equiv 1$) we deduce that problem (\ref{p1}) has maximal and minimal solutions in $[\alpha,\beta]$. However we will show that this problem does not have the extremal solutions in $[\alpha,\beta]$. \\

The family $x_{\lambda}(t)=\lambda \cos t$, $t\in I_0$, with $\lambda \in [-1,1]$, defines a set of solutions of problem (\ref{ex1}) such that $\alpha \le x_{\lambda} \le \beta$ for each $\lambda \in [-1,1]$. Notice that the zero solution is neither the least nor the greatest solution of (\ref{ex1}) in $[\alpha,\beta]$. Now let $\hat{x} \in [\alpha,\beta]$ be an arbitrary solution of problem (\ref{ex1}) and let us prove that $\hat x$ is neither the least nor the greatest solution of (\ref{ex1}) in $[\alpha,\beta]$. First, if $\hat{x}$ changes sign in $I_0$ then $\hat{x}$ cannot be a extremal solution of problem (\ref{ex1}) because it cannot be compared with the solution $x \equiv 0$. If, on the other hand, $\hat{x} \ge 0$ in $I_0$ then the differential equation yields $\hat{x}' \le 0$ a.e. on $I_0$, which implies, along with the initial condition $\hat{x}(-\frac{\pi}{2})=0$, that $\hat{x}(t)=0$ for all $t \in I_0$. Reasoning in the same way, we can prove that $\hat{x}\le 0$ in $I_0$ implies $\hat x \equiv 0$. Hence problem (\ref{ex1}) does not have extremal solutions in $[\alpha,\beta]$.
\end{example}

The previous example notwithstanding, existence of extremal solutions for problem (\ref{p1}) between given lower and upper solutions can be proven under a few more assumptions. Specifically, we have the following extremality result.

\begin{theorem}\label{main2} Consider the problem
\begin{equation}\label{p2}
\left\{
\begin{array}{ll}
x'(t)=f(t,x(t),x(\tau(t))) \ \mbox{for a.a. } t \in I_0, \\
\\
x(t)=\Lambda(x)+k(t) \ \mbox{for all $t \in I_-$.}
\end{array}
\right.
\end{equation}

If (\ref{p2}) satisfies all the conditions in Theorem \ref{main1} and, moreover, $f$ is nondecreasing with respect to  its third variable and $\Lambda$ is nondecreasing in $[\alpha,\beta]$, then problem (\ref{p2}) has the extremal solutions in $[\alpha,\beta]$.
\end{theorem}

\noindent
{\bf Proof. } Theorem \ref{main1} guarantees that problem (\ref{p2}) has a nonempty set of solutions between $\alpha$ and $\beta$.  We will show that this set of solutions is, in fact, a directed set, and then we can conclude that it has the extremal elements by virtue of \cite[Theorem 1.2]{cid}.

According to Remark \ref{rem1}, the lower solution $\alpha$ and the upper solution $\beta$ satisfy, respectively, inequalities (\ref{sub2}) and (\ref{sobre2}) and, conversely, if $\alpha$ and $\beta$ satisfy (\ref{sub2}) and (\ref{sobre2}) then they are lower and upper solutions in the sense of Definition \ref{sub-sobre}.

Let $x_1,x_2 \in [\alpha,\beta]$ be two solutions of problem (\ref{main2}). We are going to prove that there is a solution $x_3 \in [\alpha,\beta]$ such that $x_i \le x_3$ ($i=1,2$), thus showing that the set of solutions in $[\alpha,\beta]$ is upwards directed. To do so, we consider the function $\hat{x}(t)=\max\{x_1(t),x_2(t)\}, \ t \in I_0,$ which is absolutely continuous on $I_0$. For a.a. $t \in I_0$ we have either
$$
\hat{x}'(t)=f(t,\hat{x}(t),x_1(\tau(t))),$$
or
$$
\hat{x}'(t)=f(t,\hat{x}(t),x_2(\tau(t))),$$
and, since $f$ is nondecreasing with respect to its third variable, we obtain
$$
\hat{x}'(t)=f(t,\hat{x}(t),\hat{x}(\tau(t))).$$
We also have $\hat x(t) \le \Lambda(\hat x)+k(t)$ in $I_-$ because $\Lambda$ is nondecreasing,
so $\hat{x}$ is a lower solution for problem (\ref{p2}). Theorem \ref{main1} ensures now that (\ref{p2}) has at least one solution $x_3 \in [\hat{x},\beta]$.

Analogous arguments show that the set of solutions of (\ref{p2}) in $[\alpha,\beta]$ is downwards directed and, therefore, it is a directed set.
\qed

Next we show the applicability of Theorem \ref{main2}.
\begin{example} Let $L >0$ and consider the following differential equation with reflection of argument and a singularity
at $x=0$:

\begin{equation}\label{exsingular}
x'(t)=\dfrac{-t}{x(-t)} \ \mbox{for a.a. $t \in [0,L]$}, \quad x(t)=k(t)=t \cos t -3t \ \mbox{for all $t \in [-L,0]$.}
\end{equation}
In this case, the function defining the equation is $f(t,y)=\dfrac{-t}{y},$ which is nondecreasing with respect to $y$. On the other hand, functions
$$
\alpha(t)=\left\{
\begin{array}{cl}
-2t, \ &\mbox{if $t <0$}, \\
\\
-\dfrac{1}{2}t, \ &\mbox{if $0 \le t \le L$},
\end{array}
\right.
$$
and
$$
\beta(t)=\left\{
\begin{array}{cl}
-4t, \ &\mbox{if $t <0$}, \\
\\
0, \ &\mbox{if $0 \le t \le L$},
\end{array}
\right.
$$
are lower and upper solutions for problem (\ref{exsingular}). Indeed, for $t \in [-L,0]$ we have $-4t \le k(t) \le -2t$ and for a.a. $t \in I_0$ we have
$$
f(t,\alpha(-t))=-\dfrac{1}{2}=\alpha'(t), \quad f(t,\beta(-t))=-\dfrac{1}{4}<\beta'(t).$$

Finally, for a.a. $t \in I_0$ and all $y \in [\alpha(-t),\beta(-t)]$ we have
$$
|f(t,x,y)| \in \left[\dfrac{1}{4},\dfrac{1}{2}\right],$$
so problem (\ref{exsingular}) has the extremal solutions  in $[\alpha,\beta]$. Notice that $f$ admits a Carath\'eodory extension to $I_0\times \R$ outside the set
$$\{(t,y) \in I_0 \times \R \, : \, \alpha(-t) \le y \le \beta(-t)\},$$
so Theorem \ref{main2} can be applied.

In fact, we can explicitly solve problem (\ref{exsingular}) because the differential equation and the initial condition yield
$$x'(t)=\dfrac{1}{\cos t-3} \quad \mbox{for all $t \in [0,L]$, and $x(0)=0$},$$
hence problem (\ref{exsingular}) has a unique solution which  is given by
$$x(t)=\int_0^t
{\dfrac{dr}{\cos r-3}}, \ t \in [0,L].$$

\begin{figure}[H]
  \centering
    \includegraphics[width=0.4\textwidth]{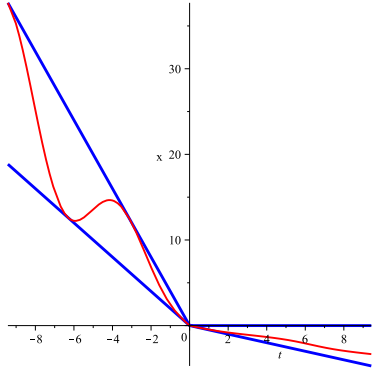}
  \caption{Solution of (\ref{exsingular}) bracketed by the lower and the upper solution.}
  \label{grafica1}
\end{figure}
\end{example}

\section{Construction of lower and upper solutions}

In general, condition $(H_1)$ is the most difficult to check among all the hypotheses in Theorem \ref{main1}. Because of this, we include in this section some sufficient conditions on the existence of linear lower and upper solutions for problem (\ref{p1}) in particular cases. \\

We begin by considering a problem of the form
\begin{equation}\label{subsol1}
\left\{
\begin{array}{ll}
x'(t)=f(x(\tau(t,x))) \ \mbox{ for a.a. } t \in I_0=[t_0,t_0+L], \\
\\
x(t)=k(t) \ \mbox{for all } t \in I_-=[t_0-r,t_0],
\end{array}
\right.
\end{equation}
where $f \in \mathcal{C}(\R)$ and $k \in \mathcal{C}(I_-)$.

\begin{proposition}\label{lowup1} Assume that $f$ is a continuous function satisfying
\begin{eqnarray}
\label{lowup11}
\lim_{y \to + \infty} f(y)=+ \infty; \\
\label{lowup12}
\lim_{y \to -\infty} f(y)=-\infty; \\
\label{lowup13}
\lim_{y \to \pm \infty} \frac{f(y)}{y} < \frac{1}{L}.
\end{eqnarray}

Then there exist $m,\overline{m} >0$ such that the functions
\begin{equation}\label{alpha1}
\alpha(t)=\left\{
\begin{array}{ll}
\varphi_*, \ &\mbox{if } t < t_0, \\
m(t_0-t)+ \varphi_*, \ &\mbox{if } t \ge t_0,
\end{array}
\right.
\end{equation}
and
\begin{equation}\label{beta1}
\beta(t)=\left\{
\begin{array}{ll}
\varphi^*, \ &\mbox{if } t < t_0, \\
\overline{m}(t-t_0)+ \varphi^*, \ &\mbox{if } t \ge t_0,
\end{array}
\right.
\end{equation}
are, respectively, a lower and an upper solution for problem (\ref{subsol1}), where
$$
\varphi_*= \min_{t \in I_-} k(t), \quad \varphi^*= \max_{t \in I_-} k(t).$$

In particular, problem (\ref{subsol1}) has maximal and minimal solutions between $\alpha$ and $\beta$, and this does not depend on the choice of $\tau$.

\end{proposition}

\noindent {\bf Proof.} Conditions (\ref{lowup12}) and (\ref{lowup13}) imply that
$$
\lim_{y \to - \infty} \frac{y-\varphi_*}{f(y)}>L,$$
so there exists $y_1 < \min\{0,\varphi_*\}$ such that
\begin{equation}\label{ine1}
0 > f(y) > \frac{y-\varphi_*}{L} \ \mbox{if } y \le y_1.
\end{equation}
On the other hand, condition (\ref{lowup11}) implies that there exists $y_2 > 0$ such that
\begin{equation}\label{ine2}
f(y) > 0 \ \mbox{if } y \ge y_2.
\end{equation}

Let $\lambda = \min\{f(y) \, : \, y_1 \le y \le y_2\}$. By condition (\ref{lowup12}) and continuity of $f$, there exists $y_3 \le y_1$ such that
\begin{equation}\label{ine3}
f(y_3)=\lambda \ \mbox{and } f(y) \ge \lambda \ \mbox{for all } y \in [y_3,y_1],
\end{equation}
and this choice of $y_3$ also provides that
\begin{equation}\label{ine4}
f(y_3) \le f(y) \ \mbox{for all } y \ge y_3,
\end{equation}
and, by virtue of (\ref{ine1}),
\begin{equation}\label{ine5}
f(y_3) > \frac{y_3-\varphi_*}{L}.
\end{equation}

Now, define $\alpha$ as in (\ref{alpha1}), with $m=\frac{\varphi_*-y_3}{L}$. Notice that $\alpha(t) \ge k(t)$ for all $t \in I_-$, $\alpha'(t)=\frac{y_3-\varphi_*}{L}$ for all $t \in I_0$ and
$$
\min_{t \in I} \alpha(t)=\alpha(t_0+L)=-mL + \varphi_*=y_3,$$
so we deduce from (\ref{ine4}) and (\ref{ine5}) that for all $t \in I_0$ we have
\begin{equation}\label{ine6}
\alpha'(t)=-m < f(y_3)=\min_{y \ge \min_{I}\alpha(t)} f(y).
\end{equation}
In the same way, we can find $\overline{y}_3 \ge \max\{0,\varphi^*\}$ such that $\beta$ defined as in (\ref{beta1}) with $\overline{m}=\frac{\varphi^*-\overline{y}_3}{L}$ satisfies that $\beta(t) \ge k(t)$ for all $t \in I_-$ and
\begin{equation}\label{ine7}
\beta'(t)=\overline{m} \ge \max_{y \le \max_I \beta(t)} f(y) \ \mbox{for all } t \in I_0.
\end{equation}

So we deduce from (\ref{ine6}) and (\ref{ine7}) that $\alpha$ and $\beta$ are lower and upper solutions for problem (\ref{subsol1}). \qed

\begin{example}\label{ex2} The function
$$
f(y)= \left\{
\begin{array}{ll}
sgn(y) \, log|y|, \ &\mbox{if } y \in (-\infty,-1) \cup (1,\infty), \\
sin(\pi y), \ &\mbox{if } y \in [-1,1],
\end{array}
\right.
$$
satisfies all the conditions in Proposition \ref{lowup1} for every compact interval $I_0$. So the corresponding problem (\ref{subsol1}) has at least one solution for any choice of $k \in \mathcal{C}(I_-)$ and $\tau \in \mathcal{C}(I,I)$.

\end{example}

We use now the ideas of Proposition \ref{lowup1} to construct lower and upper solutions for the general problem (\ref{p1}).

\begin{proposition}\label{lowup2} Let $k \in \mathcal{C}(I_0)$ and let $f:I_0 \times \R^2 \longrightarrow \R$ be a Carath\'eodory function. Assume that there exist $F_{\alpha},F_{\beta} \in \mathcal{C}(\R)$ such that for a.a. $t \in I_0$ and all $y \in \R$ we have
\begin{equation}\label{lowup21}
f(t,x,y) \ge F_{\alpha}(y) \ \mbox{for all } x \le \varphi_*
\end{equation}
and
\begin{equation}\label{lowup22}
f(t,x,y) \le F_{\beta}(y) \ \mbox{for all } x \ge \varphi^*.
\end{equation}
Moreover, assume that the next conditions involving $F_{\alpha}$ and $F_{\beta}$ hold:
\begin{eqnarray}
\label{lowup23} \lim_{y \to -\infty} F_{\alpha}(y)=-\infty, \\
\label{lowup24} F_{\alpha} \ \mbox{is bounded from below in } [0,+\infty), \\
\label{lowup25} \lim_{y \to -\infty} \frac{F_{\alpha}(y)}{y} < \frac{1}{L}, \\
\label{lowup26} \lim_{y \to +\infty} F_{\beta}(y)=+\infty, \\
\label{lowup27} F_{\beta} \ \mbox{is bounded from above in } (-\infty,0], \\
\label{lowup28} \lim_{y \to +\infty} \frac{F_{\beta}(y)}{y} < \frac{1}{L}.
\end{eqnarray}

Then there exist $m,\overline{m} \ge 0$ such that $\alpha$ and $\beta$ defined as in (\ref{alpha1})--(\ref{beta1}) are lower and upper solutions for problem (\ref{p1}), and this does not depend on the choice of $\tau$.
\end{proposition}

\noindent {\bf Proof.} Reasoning in the same way that in the proof of Proposition \ref{lowup1}, we obtain that there exists $m \ge 0$ such that $\alpha(t) \le \varphi_*$ for all $t \in I_-$ and
$$
\alpha'(t) =-m \le \min_{y \ge \min_{I} \alpha} F_{\alpha}(y) \ \mbox{for a.a. } t \in I_0.$$
As $\alpha(t) \le \varphi_*$ for all $t \in I$, we obtain by virtue of (\ref{lowup21}) that
$$
\alpha'(t) \le \min_{y \ge \min_{I} \alpha} f(t,\alpha(t),y) \ \mbox{for a.a. } t \in I_0.$$
In the same way, there exists $\overline{m} \ge 0$ such that $\beta(t) \ge \varphi^*$ for all $t \in I_-$ and
$$
\beta'(t) = \overline{m} \ge \max_{y \le \max_{I} \beta} f(t,\beta(t),y) \ \mbox{for a.a. } t \in I_0.$$

Therefore, $\alpha$ and $\beta$ are lower and upper solutions for problem (\ref{p1}). \qed

\begin{example}\label{ex3} Let $F$ the function defined in Example \ref{ex2} and consider the problem
\begin{equation}\label{example}
\left\{
\begin{array}{ll}
x'(t)=-(x+\pi)|x+\pi|^{\gamma} g(t,x) + F(\tau(t,x)) \ \mbox{for a.a. } t \in [0,L], \\
x(t)=-t \cos t \ \mbox{for all } t \in [-\pi,0],
\end{array}
\right.
\end{equation}
where $\gamma \ge 0$, $L>0$, and $g$ is a nonnegative Carath\'eodory function. \\

In this case we have $\varphi_*=-\pi$, $\varphi^*\approx 0.5611$, and the function $f(t,x,y)$ which defines the equation satisfies
$$
f(t,x,y) \ge F(y) \ \mbox{if } x \le -\pi \quad \mbox{and } f(t,x,y) \le F(y) \ \mbox{if } x \ge -\pi,$$
so in particular conditions (\ref{lowup21}) and (\ref{lowup22}) hold. As conditions (\ref{lowup23})--(\ref{lowup28}) also hold (see Example \ref{ex2}) we obtain that there exist $m,\overline{m} >0$ such that $\alpha$ and $\beta$ defined as in (\ref{alpha1})--(\ref{beta1}) are lower and upper solutions for problem (\ref{example}) for any choice of $\tau$. In particular, if there exists $\psi \in L^1(I_0)$ such that for a.a. $t \in I_0$ and all $x \in [\alpha(t),\beta(t)]$ we have $g(t,x) \le \psi(t)$, then problem (\ref{example}) has maximal and minimal solutions between $\alpha$ and $\beta$.

\end{example}

\begin{remark} Notice that the lower and upper solutions obtained both in Propositions \ref{lowup1} and \ref{lowup2} satisfy a slightly stronger condition than the one required  in Definition \ref{sub-sobre}.
\end{remark}

\end{document}